\documentclass[12pt,thmsa]{article}
\usepackage{sw20bams}

\input tcilatex
\begin{document}

\author{Steven R. Finch}
\title{Several Constants Arising in Statistical Mechanics}
\date{January 29, 1999}
\maketitle

\begin{abstract}
This is a brief survey of certain constants associated with random lattice
models, including self-avoiding walks, polyominoes, the Lenz-Ising model,
monomers and dimers, ice models, hard squares and hexagons, and percolation
models.
\end{abstract}

\section{Introduction}

Random lattice models give rise to combinatorial problems that are often
easily-stated but intractable. This article briefly presents several such
problems, each involving the numerical estimation of asymptotic growth
constants. Emphasis is given to recent developments. The author has been
collecting and writing about various mathematical constants for five years.
His large evolving website \cite{Finch} provides further discussion and
references; updates and corrections from readers are always welcome.

\section{Self-Avoiding Walks}

Let $L$ denote the $d$-dimensional cubic lattice whose sites (vertices) are
precisely all integer points in $d$-dimensional space. An $n$-step \textbf{%
self-avoiding walk} $\omega $ on $L$, beginning at the origin, is a sequence
of sites $\omega (0)$, $\omega (1)$, $\omega (2)$, ..., $\omega (n)$ with $%
\omega (0)=0$, $\left| \omega (j+1)-\omega (j)\right| =1$ for all $j$ and $%
\omega (i)\neq \omega (j)$ for all $i\neq j$.

The number of such walks $\omega $ is denoted by $c(n)$. For example, $%
c(0)=1 $, $c(1)=2d$ and $c(2)=2d(2d-1)$.

What can be said about the asymptotics of $c(n)$? It is known that 
\[
\mu =\lim\limits_{n\rightarrow \infty }c(n)^{\frac 1n} 
\]
exists and is nonzero. The value $\mu $ is called the \textbf{connective
constant} and clearly depends on $d$. The current best rigorous lower and
upper bounds for $\mu $, plus the best known estimates for $\mu $, are given
in the following table \cite{MadrasSlade}, \cite{Alm}, \cite{Noonan}, \cite
{ConwayGuttmann}, \cite{HaraSladeSokal}: 
\[
\begin{tabular}{llll}
$d$ & lower bound & best estimate & upper bound \\ 
2 & 2.62002 & 2.6381585 & 2.6939 \\ 
3 & 4.572140 & 4.683907 & 4.7476 \\ 
4 & 6.742945 & 6.7720 & 6.8179 \\ 
5 & 8.828529 & 8.8386 & 8.8602 \\ 
6 & 10.874038 & 10.8788 & 10.8886
\end{tabular}
\]
The upper bounds were computed only recently \cite{Noonan} via the
Goulden-Jackson cluster method \cite{NoonanZeilberger}. Similar techniques
can be used to estimate other constants associated with the combinatorics of
words, e.g., the asymptotics of ternary square-free words and of binary
cube-free words \cite{Finch}.

The connective constant values $\mu $ given above apply not only to the
growth of the number of self-avoiding walks, but also to the growth of
numbers of \textit{self-avoiding polygons} and of self-avoiding walks with 
\textit{prescribed endpoints} \cite{Hughes}.

For $d=2$ and $3$, there apparently is a positive constant $\gamma $ such
that

\[
\lim_{n\rightarrow \infty }\frac{c(n)}{\mu ^n\cdot n^{\gamma -1}} 
\]
exists and is nonzero \cite{MadrasSlade}, \cite{ButeraComi}, \cite{CCP}.
These constants are conjectured to be $\frac{43}{32}$ if $d=2$ and $%
1.1575... $ if $d=3$.

Another interesting object of study is the mean square displacement

\[
s(n)=E\left\{ \left| \omega (n)\right| ^2\right\} =\frac 1{c(n)}\cdot
\sum_\omega \left| \omega (n)\right| ^2 
\]
where the summation is over all $n$-step self-avoiding walks $\omega $ on $L$%
, weighted with uniform probability. Like $c(n)$, it's believed for $d=2$
and $3$ that there is a positive constant $\nu $ such that

\[
\lim_{n\rightarrow \infty }\frac{s(n)}{n^{2\nu }} 
\]
exists and is nonzero; moreover \cite{MadrasSlade}, \cite{ButeraComi}, \cite
{LMS} it's conjectured that $\nu =\frac 34$ if $d=2$ and $\nu =0.5877...$ if 
$d=3$. The \textbf{critical exponents} $\gamma $ and $\nu $ are thought to
be \textit{universal} in the sense that they are lattice-independent
(although dimension-dependent). No one has yet discovered a proof of their
existence on $L$, let alone a proof of universality.

\section{Polyminoes}

A \textbf{domino} is a pair of adjacent squares. 
\[
\text{\TEXTsymbol{<}FIGURE 1\TEXTsymbol{>}} 
\]
Generalizing, a \textbf{polyomino} or \textbf{lattice animal} of order $n$
is a connected set of $n$ adjacent squares, e.g., for $n=3$, 
\[
\text{\TEXTsymbol{<}FIGURE 2\TEXTsymbol{>}} 
\]
Define $A(n)$ to be the number of polyminoes of order $n$, where it's agreed
that two polyminoes are distinct iff they have different shapes \textit{or}
different orientations. There are different senses in which polyminoes are
defined, e.g., free versus fixed, bond versus site, simply-connected versus
not necessarily so, and others. For brevity's sake, we focus only on the
fixed, site, possibly multiply-connected case.

Redelmeier \cite{Redelmeier} computed $A(n)$ up to $n=24$ and Conway and
Guttmann \cite{ConwayGuttmann2} extended the sequence to $n=25$. It is known
that the limit 
\[
\alpha =\lim_{n\rightarrow \infty }A(n)^{\frac 1n} 
\]
exists and is nonzero. The best known bounds on $\alpha $ are $3.791\leq
\alpha \leq 4.649551$ as discussed in \cite{KlarnerRivest}, \cite{RandsWelsh}%
, \cite{WhittingtonSoteros}, \cite{Guttmann}. Improvements are possible
using the new value $A(25)$. The best known estimate, obtained via series
expansion analysis by differential approximants \cite{ConwayGuttmann2}, is $%
\alpha =4.06265...$ A more precise asymptotic expression \cite
{ConwayGuttmann2} for $A(n)$ is

\[
C\cdot n^{-1}\cdot \alpha ^n 
\]
where $C=0.316...$, but such an empirical result is far from being
rigorously proved.

\section{Lenz-Ising Model}

Let $L_N$ denote the regular $d$-dimensional cubic lattice with $N=n^d$
sites. For example, in two dimensions, $L_N$ is the $n\times n$ square
lattice with $N=n^2$. To eliminate boundary effects, $L_N$ is wrapped around
to form a $d$-dimensional torus so that, without exception, every site has $%
2d$ nearest neighbors.

Let's agree that a nonempty \textit{subgraph} of $L_N$ is connected and
contains at least one bond (edge). Suppose that several subgraphs are drawn
on $L_N$ with the property that

\begin{itemize}
\item  each bond of $L_N$ is used at most once

\item  each site of $L_N$ is used an \textit{even} number of times (possibly
zero).
\end{itemize}

\noindent Call such a configuration on $L_N$ an \textbf{even polygonal
drawing}.

\[
\text{\TEXTsymbol{<}FIGURE 3\TEXTsymbol{>}} 
\]
Note that an even polygonal drawing is the union of simple, closed,
bond-disjoint polygons but need not be connected. Other names in the
literature for these configurations include closed or Eulerian subgraphs.
Let $B(r)$ be the number of even polygonal drawings for which there are
exactly $r$ bonds. For example, when $d=2$, it follows \cite{Cipra} that $%
B(4)=N$, $B(6)=2N$ and $B(8)=\frac 12N(N+9)$ . When $d=3$, it follows that $%
B(4)=3N$, $B(6)=22N$ and $B(8)=\frac 12N(9N+375)$ . Computing $B(r)$ for
larger $r$ is quite complicated, especially considering that in $>2$
dimensions, the drawings can intertwine and be knotted!

Define the \textbf{high temperature zero magnetic field free energy} for the 
\textbf{Ising model} to be the series

\[
\beta (z)=\lim_{n\rightarrow \infty }\frac 1N\cdot \ln (1+\sum_rB(r)\cdot
z^r)=\sum_{k=0}^\infty \beta _k\cdot z^k 
\]
where $z$ is called the \textbf{activity}. If $T$ denotes temperature, then $%
\tanh ^{-1}(z)\propto 1/T$ . Coefficients $\beta _k$ of this series, as
polynomials in $d$, are given here \cite{Harris}:

\begin{center}
\begin{tabular}{ll}
$k$ & $\beta _k$ \\ 
4 & $\frac 12d(d-1)$ \\ 
6 & $\frac 13d(d-1)(8d-13)$ \\ 
8 & $\frac 14d(d-1)(108d^2-424d+425)$ \\ 
10 & $\frac 2{15}d(d-1)(2976d^3-19814d^2+44956d-34419)$%
\end{tabular}
\end{center}

The radius of convergence \cite{GuttmannEnting}, \cite{BloteLuijten}, \cite
{MunkelHeermann}, \cite{GofmanAdler}, \cite{ButeraComi} for $\beta (z)$ is 
\[
z_c=\lim_{k\rightarrow \infty }\beta _{2k}^{-\frac 1{2k}}=\left\{ 
\begin{array}{ccc}
\sqrt{2}-1=0.414213562373095049... & \text{if} & d=2 \\ 
0.218094... & \text{if} & d=3 \\ 
0.14855... & \text{if} & d=4 \\ 
0.1134... & \text{if} & d=5 \\ 
0.0920... & \text{if} & d=6 \\ 
0.0775... & \text{if} & d=7
\end{array}
\right. 
\]
which is important since knowledge of $z_c$ gives the \textit{critical
temperature} or \textit{Curie point} $T_c$ of the model. The exact two
dimensional result is a famous outcome of work by Kramers and Wannier and by
Onsager.

Here is a related problem. Suppose that several subgraphs are drawn on $L_N$
with the property that

\begin{itemize}
\item  each bond of $L_N$ is used at most once

\item  all sites of $L_N$, except two, are even

\item  the two remaining sites are odd and must lie in the same (connected)
subgraph.
\end{itemize}

\noindent Call this configuration an \textbf{odd polygonal drawing}. Note
that an odd polygonal drawing is the bond-disjoint union of an even
polygonal drawing and an (undirected) self-avoiding walk linking the two odd
sites. 
\[
\text{\TEXTsymbol{<}FIGURE 4\TEXTsymbol{>}} 
\]
Enumerating odd polygonal drawings gives rise to what physicists call the 
\textbf{high temperature zero field magnetic susceptibility } 
\[
\chi (z)=\sum_{k=0}^\infty \chi _k\cdot z^k 
\]
Coefficients $\chi _k$ of this series, as polynomials in $d$, are listed in 
\cite{GofmanAdler}. The radius of convergence of $\chi (z)$ is the same as
that for $\beta (z)$ for $d>1$. Further, when $d=2$, 
\[
\lim_{T\rightarrow T_c^{+}}\left( 1-\frac{T_c}T\right) ^{\frac 74}\cdot \chi
(z)=0.9625817322... 
\]
Wu, McCoy, Tracy, Barouch \cite{WuMcCoyTracyBarouch}, \cite
{GartenhausMcCullough} determined an exact formula for this coefficient in
terms of the Painlev\'e III function. The expression is complicated: can a
simpler formula in terms of other mathematical constants (e.g., Glaisher's
constant \cite{Finch}) be found? An exact analog of this formula for $d=3$
is also evidently not known.

\section{Monomers and Dimers}

Let $L_N$ be as before, but without wraparound. Two sites of $L_N$ are
called adjacent if the distance between them is 1. A \textbf{dimer} consists
of two adjacent sites of $L_N$ and the (non-oriented) bond connecting them.
A \textbf{dimer arrangement} is a collection of disjoint dimers on $L_N$ .
Uncovered sites are called \textbf{monomers}, so dimer arrangements are also
known as \textbf{monomer-dimer coverings}. A \textbf{dimer covering} is a
dimer arrangement whose union contains all the sites of $L_N$.

For $d=2$, let $g(n)$ denote the number of distinct monomer-dimer coverings
of $L_N$, then clearly $g(1)=1$, $g(2)=7$ and asymptotically \cite{Baxter1}, 
\cite{Heise}

\[
\kappa =\lim_{n\rightarrow \infty }g(n)^{\frac 1N}=1.940215351... 
\]
No exact expression for the constant $\kappa $ is known. Baxter's approach
for estimating $\kappa $ was based on the corner transfer matrix variational
approach. A natural way for physicists to discuss the monomer-dimer problem
is to introduce an activity $z$ for the number of monomers. The constant $%
\kappa $ then corresponds to the situation in which $z=1$. Values of $g(n)$
were recently computed \cite{Henry} up to $n=21$.

The above contrasts with the special case of dimer coverings. An exact
expression is known here for $d=2$, due to Kastelyn, Fisher and Temperley.
If $f(n)$ is the number of distinct dimer coverings of $L_N$, then $f(n)=0$
if $n$ is odd and asymptotically \cite{Percus}

\[
\lim\Sb n\rightarrow \infty  \\ n\text{ even}  \endSb f(n)^{\frac 2N}=\exp (%
\frac{2G}\pi )=1.79162281206959342... 
\]
a fascinating and unexpected occurrence of Catalan's constant $G$. For $d=3$%
, the number $h(n)$ of distinct dimer coverings \cite{Lundow} of $L_N$ is $%
h(2)=9$, $h(4)=5051532105$ and asymptotically

\[
\lim\Sb n\rightarrow \infty  \\ n\text{ even}  \endSb h(n)^{\frac 2N}=\exp
(\lambda ) 
\]
where the constant $\lambda $ is known only imprecisely. The current best
rigorous bounds \cite{Priezzhev}, \cite{Ciucu}, \cite{Schrijver} are $%
0.44007584\leq \lambda \leq 0.463107$ and the best known estimate \cite
{BeichlSullivan} is $\lambda =0.4466...$ .

\section{Ice Models}

Let $L_N$ be as before, with wraparound. An \textbf{orientation} of $L_N$ is
an assignment of a direction (arrow) to each bond of $L_N$. Assume that $d=2$
henceforth. What is the number, $\theta (n)$, of orientations of $L_N$ such
that at each site there are exactly two inward and two outward pointing
edges? (Such orientations are said to obey the \textbf{ice rule} and are
also called \textbf{Eulerian orientations}.) Here is a sample configuration:

\[
\text{\TEXTsymbol{<}FIGURE 5\TEXTsymbol{>}} 
\]
After intricate analysis, Lieb proved that \cite{Percus}, \cite{Stanley} 
\[
\lim_{n\rightarrow \infty }\theta (n)^{\frac 1N}=\left( \frac 43\right)
^{\frac 32}=1.539600717839002039... 
\]
This constant is known as the \textbf{residual entropy for square ice}.
Interestingly, $3\theta (n)$ is also the number of ways of coloring the
square faces of $L_N$ with three colors so that no two adjacent faces are
colored alike \cite{Baxter2}. Values of $\theta (n)$ were recently computed 
\cite{Henry} up to $n=13$, and relevant computational complexity issues were
discussed in \cite{MihailWinkler}.

The residual entropy $W$ for \textit{ordinary hexagonal ice} Ice-Ih and for 
\textit{cubic ice} Ice-Ic (both complicated three-dimensional lattices)
satisfy \cite{Nagle} 
\[
1.5067<W<1.5070 
\]
and are equal within the limits of Nagle's estimation error. These
configurations are not the same as the simple models mathematicians tend to
focus on. It would be interesting to see the value of $W$ for the customary $%
n\times n\times n$ cubic lattice $L_3$, either with the ice rule in effect
(two arrows point out and four arrows point in) or with Eulerian orientation
(three arrows point out and three arrows point in). No one appears to have
done this.

\section{Hard Squares and Hexagons}

Consider the set of all $n\times n$ binary matrices. What is the number $%
F(n) $ of such matrices with no pairs of adjacent 1's? Two 1's are said to
be adjacent if they lie in positions $(i,j)$ and $(i+1,j)$, or if they lie
in positions $(i,j)$ and $(i,j+1)$, for some $i$, $j$. Equivalently, $F(n)$
is the number of configurations of non-attacking Princes on an $n\times n$
chessboard, where a ''Prince'' attacks the four adjacent, non-diagonal
places. Let $N=n^2$, then 
\[
\xi =\lim_{n\rightarrow \infty }F(n)^{\frac 1N}=1.50304808247533226... 
\]
is the \textbf{hard square entropy constant} \cite{BaxterEntingTsang}, \cite
{CalkinWilf}, \cite{McKay}, \cite{Baxter3}. Essentially nothing is known
about the arithmetic character of $\xi $.

Instead of an $n\times n$ binary matrix, consider an $n\times n$ binary
array which looks like: 
\[
\left( 
\begin{array}{cccc}
a_{11} &  & a_{23} &  \\ 
& a_{22} &  & a_{34} \\ 
a_{21} &  & a_{33} &  \\ 
& a_{32} &  & a_{44} \\ 
a_{31} &  & a_{43} &  \\ 
& a_{42} &  & a_{54} \\ 
a_{41} &  & a_{53} &  \\ 
& a_{52} &  & a_{64}
\end{array}
\right) 
\]
(here $n=4$). What is the number $G(n)$ of such arrays with no pairs of
adjacent 1's? Two 1's here are said to be adjacent if they lie in positions $%
(i,j)$ and $(i+1,j)$, or in $(i,j)$ and $(i,j+1)$, or in $(i,j)$ and $%
(i+1,j+1)$, for some $i$, $j$. Equivalently, $G(n)$ is the number of
configurations of non-attacking Kings on an $n\times n$ chessboard with
regular hexagonal cells. It's surprising that the \textbf{hard hexagon
entropy constant} 
\[
\eta =\lim_{n\rightarrow \infty }G(n)^{\frac 1N}=1.395485972479302735... 
\]
is \textit{algebraic} (in fact, is solvable in radicals \cite{Joyce}) with
minimal integer polynomial \cite{Zimmermann} 
\begin{eqnarray*}
&&25937424601x^{24}+2013290651222784x^{22}+2505062311720673792x^{20}+ \\
&&797726698866658379776x^{18}+7449488310131083100160x^{16}+ \\
&&2958015038376958230528x^{14}-72405670285649161617408x^{12}+ \\
&&107155448150443388043264x^{10}-71220809441400405884928x^8 \\
&&-73347491183630103871488x^6+97143135277377575190528x^4 \\
&&-32751691810479015985152
\end{eqnarray*}
This is a consequence of Baxter's exact solution of the hard hexagon model 
\cite{Baxter2}, \cite{Andrews} via theta elliptic functions and the
Rogers-Ramanujan identities from number theory!

Just as series for the Ising model were defined using counts of even
polygonal drawings with exactly $r$ bonds, series for the hard hexagon model
can be defined using counts of non-attacking configurations of exactly $r$
Kings (and likewise for the hard square model). The radius of convergence
for the hexagon series 
\[
z_c=\frac{11+5\sqrt{5}}2=11.09016994374947424... 
\]
possesses an exact expression. No similar theoretical breakthrough has
occurred for the square model, hence the radius of convergence for the
square series 
\[
z_c=3.7962... 
\]
has no known analogous formula \cite{Baxter2}. These values are important
since they correspond to the critical activity (e.g., temperature or
density) at which a phase transition occurs in the model.

If one replaces Princes by Kings on the chessboard with square cells, then
the corresponding constant \cite{McKay} is 1.342643951124.... A related
problem, that of enumerating the \textit{maximal} configurations of $\frac
N4 $ nonattacking Kings, was discussed in \cite{Wilf}, \cite{Larsen}.

\section{Percolation Models}

Let $M$ be a random $n\times n$ binary matrix satisfying

\begin{itemize}
\item  $m_{ij}=1$ with probability $p$, $0$ with probability $1-p$ for each $%
i,j$

\item  $m_{ij}$ and $m_{kl}$ are independent for all $(i,j)\neq (k,l).$
\end{itemize}

\noindent An $\mathbf{s}$\textbf{-cluster} is an isolated grouping of $s$
adjacent 1's in $M$, where adjacency means horizontal or vertical neighbors
(not diagonal). For example, the $4\times 4$ matrix 
\[
M=\left( 
\begin{array}{cccc}
1 & 0 & 1 & 1 \\ 
1 & 1 & 0 & 0 \\ 
0 & 1 & 0 & 1 \\ 
1 & 0 & 0 & 1
\end{array}
\right) 
\]
has one 1-cluster, two 2-clusters and one 4-cluster. The total number of
clusters is 4 in this case. For arbitrary $n$, the total cluster count is a
random variable with normalized expected value 
\[
K_n=\frac{E(C_n)}{n^2} 
\]
The limit $K_S(p)$ of $K_n$ exists as $n$ approaches infinity, and $K_S(p)$
is called the \textbf{mean cluster count per site} or the \textbf{mean
cluster density} for the \textbf{site percolation model}. It's known that $%
K_S(p)$ is twice continuously differentiable on [0,1]; further, $K_S(p)$ is
analytic on [0,1] except possibly at one point $p=p_c$. Monte Carlo
simulation and numerical Pad\'e approximants can be used to compute $K_S(p)$%
. For example \cite{ZiffFinchAdamchik}, it's known that $K_S(\frac
12)=0.065770....$

Instead of an $n\times n$ binary matrix $M$, consider a binary array $A$ of $%
2n(n-1)$ entries which looks like 
\[
A=\left( 
\begin{array}{ccccccc}
& a_{12} &  & a_{14} &  & a_{16} &  \\ 
a_{11} &  & a_{13} &  & a_{15} &  & a_{17} \\ 
& a_{22} &  & a_{24} &  & a_{26} &  \\ 
a_{21} &  & a_{23} &  & a_{25} &  & a_{27} \\ 
& a_{32} &  & a_{34} &  & a_{36} &  \\ 
a_{31} &  & a_{33} &  & a_{35} &  & a_{37} \\ 
& a_{42} &  & a_{44} &  & a_{46} & 
\end{array}
\right) 
\]
(here $n=4$). One should associate $a_{ij}$ not with a site of the $n\times
n $ square lattice (as one does for $m_{ij}$) but with a bond. An $s$%
-cluster here is an isolated, connected subgraph of the graph of all bonds
associated with 1's. For example, the array 
\[
A=\left( 
\begin{array}{ccccccc}
& 1 &  & 0 &  & 0 &  \\ 
1 &  & 0 &  & 0 &  & 0 \\ 
& 0 &  & 1 &  & 0 &  \\ 
0 &  & 1 &  & 1 &  & 0 \\ 
& 0 &  & 1 &  & 0 &  \\ 
1 &  & 0 &  & 0 &  & 0 \\ 
& 0 &  & 0 &  & 0 & 
\end{array}
\right) 
\]
has one 1-cluster, one 2-cluster and one 4-cluster. For \textbf{bond
percolation models} such as this, one often includes 0-clusters in the total
count as well, that is, isolated sites with no attached 1's bonds. In this
case there are seven 0-clusters, hence the total number of clusters $C_4$ is
10. The limiting mean cluster density 
\[
K_B(p)=\lim_{n\rightarrow \infty }K_n=\lim_{n\rightarrow \infty }\frac{E(C_n)%
}{n^2} 
\]
exists as $n$ approaches infinity and similar smoothness properties hold.
Remarkably, however, an exact integral expression exists at $p=\frac 12$ for
the limiting mean cluster density \cite{TemperleyLieb}, \cite{Essam} 
\[
K_B(\frac 12)=\left. -\frac 18\cot (y)\cdot \frac d{dy}\left\{ \frac 1y\cdot
\dint_{-\infty }^\infty \limfunc{sech}\left( \frac{\pi x}{2y}\right) \cdot
\ln \left( \frac{\cosh \left( x\right) -\cos (2y)}{\cosh \left( x\right) -1}%
\right) dx\right\} \right| _{y=\frac \pi 3} 
\]
which Adamchik \cite{ZiffFinchAdamchik} recently simplified to 
\[
K_B(\frac 12)=\frac{3\sqrt{3}-5}2=0.09807621135331594... 
\]
An analogous expression for the limiting variance of cluster density was
computed in \cite{ZiffFinchAdamchik}. In the same way, the bond percolation
model on the \textit{triangular} lattice gives a known limiting mean cluster
density at a specific value of $p$ (discussed below), but the relevant
variance is not known here.

Let's turn attention away from mean cluster density $K(p)$ and instead
toward \textbf{mean cluster size} $S(p)$. In the examples given earlier, $%
S_4=(1+2+2+4)/4=9/4$ for the site case, $S_4=(1+2+4)/3=7/3$ for the bond
case, and $S(p)$ is the limiting value of $E(S_n)$ as $n$ approaches
infinity. The \textbf{critical probability} or \textbf{percolation threshold}
$p_c$ is defined to be \cite{StaufferAharony}, \cite{Hughes} 
\[
p_c=\inf\Sb 0<p<1  \\ S(p)=\infty  \endSb p 
\]
that is, the concentration at which an infinite cluster appears in the
infinite lattice.

For site percolation on the square lattice, there are rigorous bounds \cite
{Wierman}, \cite{vanderBergErmakov} 
\[
0.556<p_c<0.679492 
\]
and an estimate \cite{Ziff} $p_c=0.5927460....$ based on extensive
simulation. Ziff \cite{ZiffFinchAdamchik} has additionally calculated that $%
K_S(p_c)=0.0275981...$ via simulation.

In contrast, for bond percolation on the square and triangular lattices,
there are exact results due to Sykes and Essam. Keston \cite{Hughes} proved
that $p_c=\frac 12$ on the square lattice, corresponding to the expression $%
K_B(\frac 12)$ above. On the triangular lattice, Wierman \cite{Hughes}
proved that 
\[
p_c=2\cdot \sin \left( \frac \pi {18}\right) =0.347296355333860698...
\]
corresponding to another exact expression \cite{BaxterTemperleyAshley}, \cite
{ZiffFinchAdamchik} 
\begin{eqnarray*}
K_B(p_c) &=&\frac{35}4-\frac 32\cdot \csc \left( \frac \pi {18}\right) =%
\frac{23}4-\frac 32\cdot \left\{ \sqrt[3]{4\cdot \left( 1+i\sqrt{3}\right) }+%
\sqrt[3]{4\cdot \left( 1-i\sqrt{3}\right) }\right\}  \\
\  &=&0.1118442752845497...
\end{eqnarray*}
Similar results apply for the hexagonal lattice by duality.

Much energy has been placed into the computation of universal exponents for
random lattice models in $d$-dimensional space (akin to what we briefly
mentioned for self-avoiding walks). The existence of such exponents is
hypothesized on the basis of both theoretical procedures (finite size
scaling and renormalization group methods) and experimental data. We shall
not attempt to discuss this important subject, but instead refer readers to 
\cite{MadrasSlade}, \cite{StaufferAharony}, \cite{Hughes}.

\begin{description}
\item  Steven R. Finch

\item  MathSoft Inc., 101 Main Street

\item  Cambridge, MA, USA 02142

\item  \textit{sfinch@mathsoft.com}
\end{description}


\begin{thebibliography}{99}
\bibitem{Alm}  Alm, S. E., Upper bounds for the connective constant of
self-avoiding walks, \textit{Combin. Probab. Comput.} 2 (1993) 115-136.

\bibitem{Andrews}  Andrews, G. E., The reasonable and unreasonable
effectiveness of number theory in statistical mechanics, \textit{Proc. Symp.
Applied Math.} 46, ed. S. A. Burr (Orono conf., 1991) Amer. Math. Soc.,
1992, 21-34.

\bibitem{Baxter1}  Baxter, R. J., Dimers on a rectangular lattice, \textit{%
J. Math. Physics} 9 (1968) 650-654.

\bibitem{Baxter2}  Baxter, R. J., \textit{Exactly Solved Models in
Statistical Mechanics}, Academic Press 1982.

\bibitem{Baxter3}  Baxter, R. J., Planar lattice gases with
nearest-neighbour exclusion, \textit{Annals Combin.} 3 (1999).

\bibitem{BaxterEntingTsang}  Baxter, R. J., Enting, I. G., Tsang, S. K.,
Hard-square lattice gas, \textit{J. Stat. Phys.} 22 (1980) 465-489.

\bibitem{BaxterTemperleyAshley}  Baxter, R. J., Temperley, H. N. V., Ashley,
S. E., Triangular Potts model at its transition temperature, and related
models, \textit{Proc. Royal Soc. London A} 358 (1978) 535-559.

\bibitem{BeichlSullivan}  Beichl, I., Sullivan, F., Approximating the
permanent via importance sampling with application to the dimer covering
problem, \textit{J. Comput. Phys.}, submitted.

\bibitem{BloteLuijten}  Bl\"ote, H. W. J., Luijten, E., Heringa, J. R.,
Ising universality in three dimensions: a Monte Carlo study, \textit{J.
Phys. Math. A} 28 (1995) 6289-6313.

\bibitem{ButeraComi}  Butera, P., Comi, M., N-vector spin models on the
simple-cubic and the body-centered-cubic lattices: a study of the critical
behavior of the susceptibility and of the correlation length by
high-temperature series extended to order 21, \textit{Phys. Rev. B} 56
(1997) 8212-8240.

\bibitem{CalkinWilf}  Calkin, N. J., Wilf, H. S., The number of independent
sets in a grid graph, \textit{SIAM J. Discrete Math. }11 (1998) 54 - 60.

\bibitem{CCP}  Caracciolo, S., Causo, M. S., Pelissetto, A., Monte Carlo
results for three-dimensional self-avoiding walks, \textit{Nucl. Phys. Proc.
Suppl.} 63 (1998) 652-654.

\bibitem{Cipra}  Cipra, B. A., An introduction to the Ising model, \textit{%
Amer. Math. Monthly} 94 (1987) 937-959.

\bibitem{Ciucu}  Ciucu, M., An improved upper bound for the three
dimensional dimer problem, \textit{Duke Math. J., }to appear.

\bibitem{ConwayGuttmann}  Conway, A. R., Guttman, A. J., Lower bound on the
connective constant for square lattice self-avoiding walks, \textit{J. Phys.
A.} 26 (1993) 3719-3724.

\bibitem{ConwayGuttmann2}  Conway, A. R., Guttmann, A. J., On
two-dimensional percolation, \textit{J. Phys. A} 28 (1995) 891-904.

\bibitem{Essam}  Essam, J. W., Percolation and cluster size, \textit{Phase
Transitions and Critical Phenomena}, vol. II, ed. C. Domb and M. S. Green,
Academic Press 1972, 197-270.

\bibitem{Finch}  Finch, S. R., \textit{Favorite Mathematical Constants},
MathSoft Inc., website URL
http://www.mathsoft.com/asolve/constant/constant.html, 1998.

\bibitem{GartenhausMcCullough}  Gartenhaus, S, McCullough, W. S., Higher
order corrections for the quadratic Ising lattice susceptibility at
criticality, \textit{Phys. Rev. B} 38 (1988) 11688-11703.

\bibitem{GofmanAdler}  Gofman, M., Adler, J., Aharony, A., Harris, A. B.,
Stauffer, D., Series and Monte Carlo study of high-dimensional Ising models, 
\textit{J. Stat. Phys.} 71 (1993) 1221-1230.

\bibitem{Guttmann}  Guttmann, A. G., On the number of lattice animals
embeddable in the square lattice, \textit{J. Phys. A} 15 (1982) 1987-1990.

\bibitem{GuttmannEnting}  Guttmann, A. J., Enting, I. G., The
high-temperature specific heat exponent of the 3D Ising model, \textit{J.
Phys. A} 27 (1994) 8007-8010.

\bibitem{HaraSladeSokal}  Hara, T., Slade, G., Sokal, A. D., New lower
bounds on the self-avoiding-walk connective constant, \textit{J. Stat. Phys.}
72 (1993) 479-517; erratum, 78 (1995) 1187-1188.

\bibitem{Harris}  Harris, A. B., Meir, Y., Recursive enumeration of clusters
in general dimension on hypercubic lattices, \textit{Phys. Rev. A} 36 (1987)
1840-1848.

\bibitem{Heise}  Heise, M., Upper and lower bounds for the partition
function of lattice models, \textit{Physica A} 157 (1989) 983-999.

\bibitem{Henry}  Henry, J. J., private communications (1997-1998).

\bibitem{Hughes}  Hughes, B. D., \textit{Random Walks and Random Environments%
}, vols. 1 and 2, Oxford, 1996.

\bibitem{Joyce}  Joyce, G. S., On the hard hexagon model and the theory of
modular functions, \textit{Phil. Trans. Royal Soc. London A} 325 (1988)
643-702.

\bibitem{KlarnerRivest}  Klarner, D. A., Rivest, R. L., A procedure for
improving the upper bound for the number of n-ominoes, \textit{Canad. J.
Math.} 25 (1973) 585-602.

\bibitem{Larsen}  Larsen, M., The problem of kings, \textit{Elec. J. Comb.}
2 (1995).

\bibitem{LMS}  Li, B., Madras, N., Sokal, A. D., Critical exponents,
hyperscaling and universal amplitude ratios for two- and three-dimensional
self-avoiding walks, \textit{J. Stat. Phys.} 80 (1995) 661-754.

\bibitem{Lundow}  Lundow, Per H\aa kan, Computation of matching polynomials
and the number of 1-factors in polygraphs, Ume\aa\ University Math. Dept.
preprint 12-1996 (1996).

\bibitem{MadrasSlade}  Madras, N., Slade, G., \textit{The Self-Avoiding Walk}%
, Birkh\"auser, 1993.

\bibitem{McKay}  McKay, B. D., private communication (1996).

\bibitem{MihailWinkler}  Mihail, M., Winkler, P., On the number of Eulerian
orientations of a graph, \textit{Proc. Third Annual ACM-SIAM Symposium on
Discrete Algorithms}, Orlando FL., 1992, pp. 138-145; also appears in 
\textit{Algorithmica} 16 (1996) 402-414.

\bibitem{MunkelHeermann}  M\"unkel, C., Heermann, D. W., Adler, J., Gofman,
M., Stauffer, D., The dynamical critical exponent of the two-, three- and
five-dimensional kinetic Ising model, \textit{Physica A} 193 (1993) 540-552.

\bibitem{Noonan}  Noonan, J., New upper bounds for the connective constants
of self-avoiding walks, \textit{J. Stat. Phys.} 91 (1998) 871-888.

\bibitem{NoonanZeilberger}  Noonan, J., Zeilberger, D., The Goulden-Jackson
cluster method: extensions, applications and implementations, \textit{J.
Difference Eq. Appl.}, to appear.

\bibitem{Nagle}  Nagle, J. F., Lattice statistics of hydrogen bonded
crystals: I. The residual entropy of ice, \textit{J. Math. Phys.} 7 (1966)
1484-1491.

\bibitem{Percus}  Percus, J. K., \textit{Combinatorial Methods},
Springer-Verlag 1971.

\bibitem{Priezzhev}  Priezzhev, V. B., The statistics of dimers on a
three-dimensional lattice, II. An improved lower bound, \textit{J. Stat. Phy.%
} 26 (1981) 829-837.

\bibitem{RandsWelsh}  Rands, B. M. I., Welsh, D. J. A., Animals, trees and
renewal sequences, \textit{IMA J. Appl. Math.} 27 (1981) 1-17.

\bibitem{Redelmeier}  Redelmeier, D. H., Counting polyominoes: Yet another
attack, \textit{Discrete Math.} 36 (1981) 191-203.

\bibitem{Schrijver}  Schrijver, A., Counting 1-factors in regular bipartite
graphs, \textit{J. Combin. Theory B} 72 (1998) 122-135; also \textit{MR}
82a:15004.

\bibitem{Stanley}  Stanley, R. P., \textit{Enumerative Combinatorics}, vol.
1, Cambridge Univ. Press, 1997.

\bibitem{StaufferAharony}  Stauffer, D., Aharony, A., \textit{Introduction
to Percolation Theory}, 2nd ed., Taylor and Francis, 1992.

\bibitem{TemperleyLieb}  Temperley, H. N. V., Lieb, E. H., Relations between
the 'percolation' and 'colouring' problem and other graph-theoretical
problems associated with regular planar lattices; some exact results for the
'percolation' problem, \textit{Proc. Royal Soc. London A} 322 (1971) 251-280.

\bibitem{vanderBergErmakov}  van den Berg, J., Ermakov, A., A new lower
bound for the critical probability of site percolation on the square
lattice, \textit{Random Structures and Algorithms} 8 (1996) 199-212.

\bibitem{WhittingtonSoteros}  Whittington, S. G., Soteros, C. E., Lattice
animals: Rigorous results and wild guesses, in \textit{Disorder in Physical
Systems: A Volume in Honour of J. M. Hammersley}, ed. G. R. Grimmett and D.
J. A. Welsh, Oxford, 1990.

\bibitem{Wierman}  Wierman, J. C., Substitution method critical probability
bounds for the square lattice site percolation model, \textit{Combin.
Probab. Comput.} 4 (1995) 181-188.

\bibitem{Wilf}  Wilf, H. S., The problem of kings, \textit{Elec. J. Comb.} 2
(1995).

\bibitem{WuMcCoyTracyBarouch}  Wu, T. T., McCoy, B. M., Tracy, C. A.,
Barouch, E., Spin-spin correlation functions for the two-dimensional Ising
model: exact theory in the scaling region, \textit{Phys. Rev. B} 13 (1976)
316-374.

\bibitem{Ziff}  Ziff, R. M., Spanning probability in 2D percolation, \textit{%
Phys. Rev. Letters} 69 (1992) 2670-2673.

\bibitem{ZiffFinchAdamchik}  Ziff, R.M., Finch, S. R., Adamchik, V.,
Universality of finite-size corrections to the number of critical
percolation clusters, \textit{Phys. Rev. Lett.} 79 (1997) 3447-3450.

\bibitem{Zimmermann}  Zimmermann, P., private communication (1996).
\end{thebibliography}
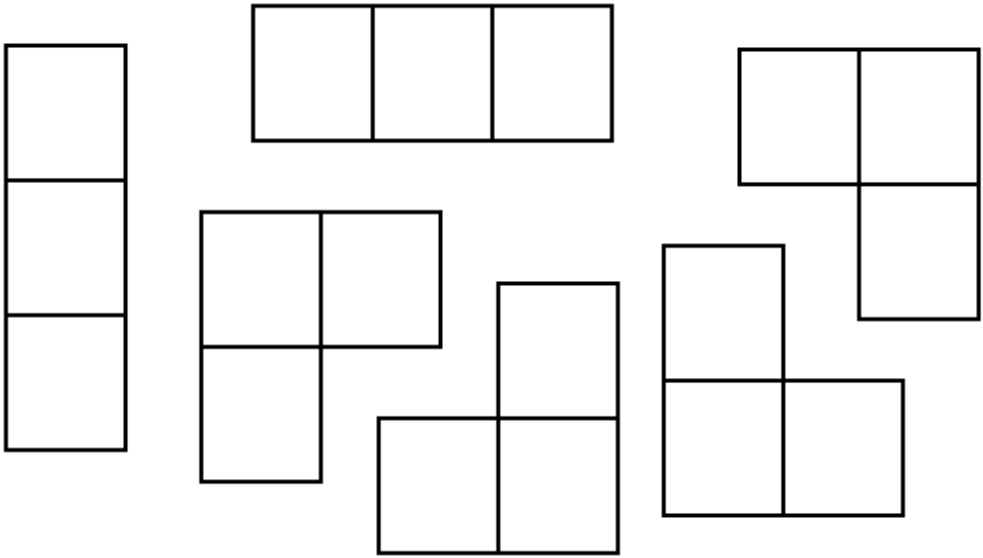
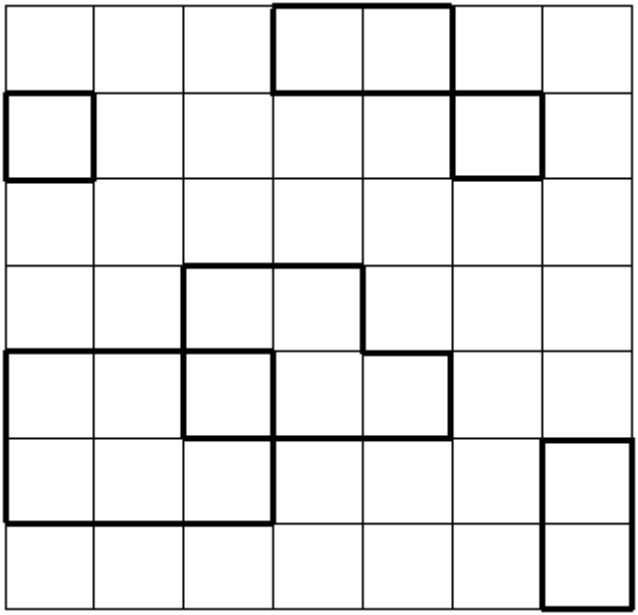
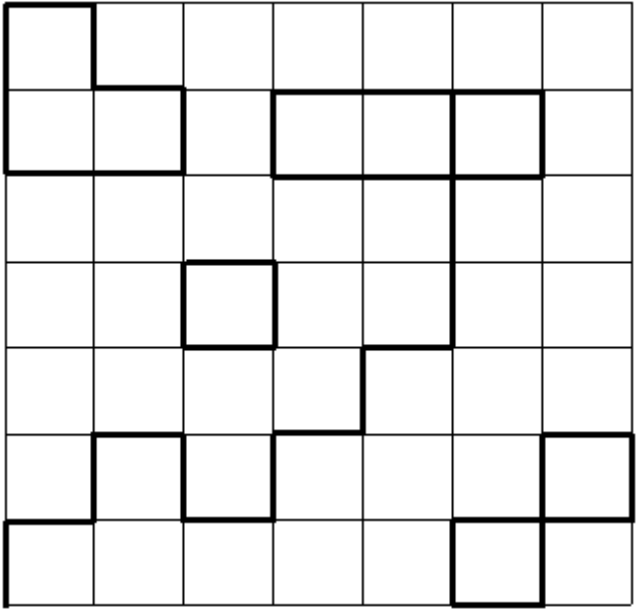
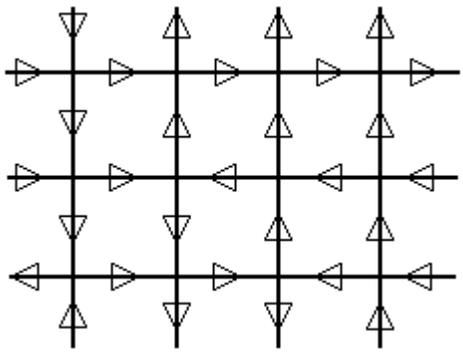
\end{document}